\documentclass[letterpaper, 9pt, conference]{ieeeconf}  

\IEEEoverridecommandlockouts                              
\usepackage{graphicx}
\usepackage{amssymb}
\usepackage{amsmath}
\usepackage{graphics} 
\usepackage{epsfig}   

\usepackage{times,verbatim} 
\usepackage{subfigure}
\usepackage{epstopdf} 
\usepackage{algorithm2e}
\usepackage[english]{babel}
\usepackage{color}

\mathchardef\Gamma="7100
\mathchardef\Delta="7101
\mathchardef\Theta="7102
\mathchardef\Lambda="7103
\mathchardef\Xi="7104
\mathchardef\Pi="7105
\mathchardef\Sigma="7106
\mathchardef\Upsilon="7107
\mathchardef\Phi="7108
\mathchardef\Psi="7109
\mathchardef\Omega="710A

\newtheorem{assumption}{Assumption}

\newtheorem{definition}{Definition}
\newtheorem{problem}{Problem}
\newtheorem{theorem}{Theorem}
\newtheorem{remark}{Remark}
\newtheorem{example}{Example}

\def\R{\mathbb{R}}

\def\fine{\hfill\hbox{$\diamond$}\medskip}
\def\fineTh{\hfill\hbox{$\blacksquare$}\medskip}
\def\nav#1{\noalign{\vskip#1pt}}
\def\crr#1{\cr\nav{#1}}

\def\norma#1{{\left\|#1\right\|}}
\def\cal{\mathcal}
\def\eqref#1{(\ref{#1})}
\def\proof{{\noindent{\itshape Proof: }}}

\def\reference{{\rm ref}}
\def\d#1{\frac{{\rm d}_+#1}{{\rm d}t}}

\def\dd#1{\frac{{\rm d}^2_+#1}{{\rm d}t^2}}
\def\dn#1{\frac{{\rm d}^n_+#1}{{\rm d}t^n}}
\long\def\Omesso#1\EndOmesso{}

\catcode`@=11
\def\eqalign#1{\null\,\vcenter{\openup\jot\m@th\ialign{\strut\hfil$\displaystyle{##}$&$\displaystyle{{}##}$\hfil\crcr#1\crcr}}\,}
\def\eqaligntwo#1{\null\,\vcenter{\openup\jot\m@th\ialign{\strut\hfil$\displaystyle{##}$&$\displaystyle{##}$\hskip.125cm&\hskip.125cm$\hfil\displaystyle{##}$&$\displaystyle{{}##}$\hfil\crcr#1\crcr}}\,}
\def\eqalignthree#1{\null\,\vcenter{\openup\jot\m@th\ialign{\strut\hfil$\displaystyle{##}$&$\displaystyle{##}$\hskip.25cm&\hskip.25cm$\hfil\displaystyle{##}$&$\displaystyle{{}##}$\hskip.25cm&\hskip.25cm$\hfil\displaystyle{##}$&$\displaystyle{{}##}$\hfil\crcr#1\crcr}}\,}
\def\eqalignno#1{\displ@y \centering\halign to\displaywidth{\hfil$\@lign\displaystyle{##}$&$\@lign\displaystyle{{}##}$\hfil\centering&\llap{$\@lign##$}\crcr#1\crcr}}
\def\tabletwo#1{\null\,\vcenter{\openup\jot\m@th\ialign{\strut\hfil$\displaystyle{##}$&$\displaystyle{##}$\hskip.125cm&\hskip.125cm$\hfil\displaystyle{##}$&$\displaystyle{{}##}$\hfil\crcr#1\crcr}}\,}
\def\matrixp#1{\null\,\vcenter{\normalbaselines\m@th
    \ialign{\hfil$##$\hfil&&$\>\>$\hfil$##$\hfil\crcr
      \mathstrut\crcr\noalign{\kern-\baselineskip}
      #1\crcr\mathstrut\crcr\noalign{\kern-\baselineskip}}}\,}
\def\pmatrix#1{\left(\matrixp{#1}\right)}

\def\cases#1{\left\{\,\vcenter{\normalbaselines\m@th
    \ialign{$##\hfil$&\quad##\hfil\crcr#1\crcr}}\right.}
\catcode`@=12

\catcode`\ˆ=\active\defˆ{\`a}
\catcode`\=\active\def{\`e}
\catcode`\˜=\active\def˜{\`o}
\catcode`\=\active\def{\`u}
\catcode`\"=\active\def"{\`\i}
\catcode`\‡=\active\def‡{\'a}
\catcode`\Ž=\active\defŽ{\'e}
\catcode`\—=\active\def—{\'o}
\catcode`\œ=\active\defœ{\'u}
\catcode`\'=\active\def'{\'\i}

\catcode`\é=\active\defé{\`E}
\catcode`\ñ=\active\defñ{\`O}
\catcode`\Ë=\active\defË{\`I}
\catcode`\ô=\active\defô{\`U}
\catcode`\Ë=\active\defË{\`A}
\catcode`\ƒ=\active\defƒ{\'E}
\catcode`\î=\active\defî{\'O}
\catcode`\ê=\active\defê{\'I}
\catcode`\ò=\active\defò{\'U}
\catcode`\ç=\active\defç{\'A}

\begin{document}

\title{\LARGE\bf Digital Self Triggered Robust Control\\of Nonlinear Systems\thanks{The research leading to these results has received funding from the European Union Seventh Framework Programme [FP7/2007-2013] under grant agreement n°257462 HYCON2 Network of excellence.}}

\author{M.D.~Di Benedetto, S.~Di Gennaro, and A.~D'Innocenzo
\thanks{\baselineskip8pt The authors are with the Department of Electrical and Information Engineering, and Center of Excellence DEWS, University of L'Aquila, Via G.~Gronchi~18, 67100, L'Aquila, Italy. E.mail: {\tt \{maria- domenica.dibenedetto,\,stefano.digennaro,\! alessandro. dinnocenzo\}@univaq.it}.}
}

\maketitle
\thispagestyle{empty}
\pagestyle{empty}

\begin{abstract}
In this paper we develop novel results on self triggering control of nonlinear systems, subject to perturbations and actuation delays. First, considering an unperturbed nonlinear system with bounded actuation delays, we provide conditions that guarantee the existence of a self triggering control strategy stabilizing the closed--loop system. Then, considering parameter uncertainties, disturbances, and bounded actuation delays, we provide conditions guaranteeing the existence of a self triggering strategy, that keeps the state arbitrarily close to the equilibrium point. In both cases, we provide a methodology for the computation of the next execution time. We show on an example the relevant benefits obtained with this approach, in terms of energy consumption, with respect to control algorithms based on a constant sampling, with a sensible reduction of the average sampling time.
\end{abstract}


\section{Introduction}\label{sec:intro}

Wireless networked control systems are spatially distributed control systems where the communication between sensors, actuators, and computational units is supported by a shared wireless communication network~\cite{Karl 2005}. The use of wireless networked control systems in industrial automation results in flexible architectures and generally reduces installation, debugging, diagnostic and maintenance costs with respect to wired networks. The main motivation for studying such systems is the emerging use of wireless technologies in control systems, see e.g. \cite{AkyildizAHN2004, SongWirelesshart2008} and references therein. Although wireless networks offer many advantages, communication nodes generally consist of battery powered devices. For this reason, when designing a control scheme closed on a wireless sensor network, it is fundamental to adopt power aware control algorithms to reduce power consumption. In other applications, the energy is obtained from the environment, with a scavenging system, since it is not possible to provide these wireless sensors with batteries. An example of such a case is notably given by the so--called intelligent (or smart) tires, equipped with sensors embedded in the tread, giving information on pressure, road--tire friction, etc.~\cite{Palmer 2002, Matsuzakia 2005, Apollo 2001}. In this case the sensors are supplied by the energy provided by the tire motion. It is clear that it is fundamental to trigger wireless transmission only when necessary, to prevent energy shortage and, possibly, to reduce the probability of information packet losses during the transmission.

With the aim of addressing the above issues in the controller design phase, self triggered control strategies have been introduced in \cite{Velasco 2003}, where a heuristic rule is provided to self--trigger the next execution time of a control task on the basis of the last measurement of the state. In~\cite{Lemmon 2007, Lemmon 2009}, a robust self triggered strategy is proposed, which guarantees that the ${\cal L}_2$ gain of a linear time invariant system is kept under a given threshold. In~\cite{Mazo 2008} a self triggering strategy distributed over a wireless sensor network is proposed for linear time invariant plants.

In \cite{Tabuada 2007} sufficient conditions for the existence of a stabilizing event--triggered control strategy are given for non-linear systems. In \cite{Anta 2010} the authors propose a self--triggered emulation of the event--triggered control strategy proposed in \cite{Tabuada 2007}. In particular a methodology for the computation of the next execution time as a function of the last sample is presented, under a homogeneity condition.

We extend the previous results in two directions. First, with regard to {\sl asymptotic stability} under digital self triggered control, we propose a methodology for the computation of the next execution time by replacing the homogeneity assumption with the requirement that the nonlinear differential equations and the control law are $C^\ell$ functions, with $\ell$ sufficiently large. The digital control, ensuring asymptotic stability, have been derived for unperturbed nonlinear systems affected by bounded actuation delays, under the necessary condition of existence of a continuous stabilizing control. With respect to previous results~\cite{Anta 2010}, the approach proposed in this work allows computing the stabilizing execution time sequence for a broader class of systems. In a recent work in progress \cite{TabuadaSubmitted2010} a new technique is developed to compute the next execution time for smooth systems by exploiting the concept of isochronous manifolds. Our technique is based on polynomial approximations of Lyapunov functions, and therefore differs from the one developed in \cite{TabuadaSubmitted2010}.

Second, we consider non--linear systems perturbed by norm--bounded parameter uncertainties and disturbances, and affected by bounded actuation delays. We prove that, under weaker conditions than those used in \cite{Tabuada 2007}, a self triggering strategy exists keeping the state in a "safe set" arbitrarily close to the equilibrium point. To the best of the authors' knowledge, this is the first work that provides results on self-triggering control for non-linear systems with uncertainties, disturbances and actuation delays.

Finally, we show on a significant example that the results obtained introduce strong benefits in terms of energy consumption, with respect to digital controls based on a constant sampling time, by reducing the average sampling time.

The paper is organized as follows. In Section~\ref{sec:Problem Formulation}, we illustrate the mathematical model and the problem formulation. In Section~\ref{sec:selfTrigStabCnt}, we derive results for the asymptotic stability of unperturbed systems, while in Section~\ref{sec:selfTrigSafetyCnt} we consider the safety problem of perturbed systems. In Section~\ref{sec:simulations}, we apply the obtained results to an example.


\section{Problem Formulation}\label{sec:Problem Formulation}

Consider a generic nonlinear system
\begin{equation}
\dot x=f(x,u,\mu,d)          \label{perturbed system}
\end{equation}
where $x\in{\cal D}_x\subset\R^n$, ${\cal D}_x$ a domain containing the origin, $u\in{\cal D}_u\subset\R^p$, $\mu$ is a parameters uncertainty vector varying in a compact set ${\cal D}_\mu\subset\R^r$, with $0\in {\cal D}_{\mu}$, $d$ is an external bounded disturbances vector taking values in a compact set ${\cal D}_d \subset\R^s$, with $0\in {\cal D}_d$. We define the nominal system associated to the perturbed system~\eqref{perturbed system} by
\begin{equation}
\dot x=f_0(x,u) \doteq f(x,u,0,0).               \label{nominal system}
\end{equation}
Given a state feedback control law $\kappa\colon{\cal D}_x\to{\cal D}_u$, the closed loop perturbed system is
\begin{equation}
\dot x=f(x,\kappa(x),\mu,d),          \label{perturbed controlled nonlinear system}
\end{equation}
and the closed loop nominal system is
\begin{equation}
\dot x=f_0(x,\kappa(x)).          \label{nominal controlled nonlinear system}
\end{equation}
We will denote by $x(t)$, $t\ge t_0$, the solution of the closed loop system~\eqref{perturbed controlled nonlinear system} (or \eqref{nominal controlled nonlinear system}, according to the context), with initial condition $x_0=x(t_0)$. It is well--known that if the origin of system~\eqref{nominal controlled nonlinear system} is locally asymptotical stable for a certain feedback $\kappa$, and if $f_0(x,\kappa(x))\in C^\ell({\cal D}_x)$, $\ell>1$ integer, then there exists a Lyapunov function $V(x)$ of class $C^1({\cal D}_x)$ such that
\begin{equation}
\eqalign{
\alpha_1(\|x\|)\le V(x)&\le\alpha_2(\|x\|) \cr
\frac{\partial V(x)}{\partial x}f_0(x,\kappa(x))&\le-\alpha_3(\|x\|) \cr
\bigg\|\frac{\partial V(x)}{\partial x}\bigg\|&\le\alpha_4(\|x\|)  \cr}   \label{Lyapunov candidate}
\end{equation}
with $\alpha_1,\alpha_2,\alpha_3,\alpha_4\in {\cal K}$~\cite{Yoshizawa 1955}, \cite{Kurzweil 1956}, \cite{Khalil 2002}.


In the following definition, an invariant property is used to define that a system is safe with respect to a given subset of the state space.

\smallskip
\begin{definition}\label{Safe}
Given a state feedback control law $\kappa$, system~\eqref{perturbed controlled nonlinear system} is safe with respect to the set ${\cal S}\subseteq {\cal D}_x$ for the time interval ${\cal T} \subseteq \R^+$, if $x(t)\in {\cal S}$, $\forall t\in {\cal T}$.\fine
\end{definition}

The feedback control signal $u(t)=\kappa(x(t))$ requires continuous measurements of the state of the system. In this paper, we address the stability and safety problems defined below, when the measurements are performed at sampling instants $t_k$, defining a sequence ${\cal I}=\{t_k\}_{k\ge0}$, and the applied control signal is
$$
u_{{\cal I}}(t)=
\cases{  0             & $\forall t \in [t_0, t_0 + \Delta_0)$ \crr{4}
        \kappa(x(t_k)) & $\forall t \in [t_k + \Delta_k, t_{k+1} + \Delta_{k+1}),\ k \ge 0$}
$$
where $\{\Delta_k\}_{k\ge 0}$ is the sequence of actuation delays, due to the transmission time from the sensor to the controller, the computation time, and the transmission time from the controller to the actuator. We assume that $\Delta_k\in[0,t_{k+1}-t_k)$, $\forall k\ge 0$, which is a natural requirement in practice.

\smallskip
\begin{problem} {\sl (Stability problem)}\label{Stability Problem}
Given the nominal system~\eqref{nominal system}, and a stabilizing state feedback control law $\kappa$, determine
\begin{enumerate}
\item[1.] A minimum sampling time $\tau_{\min} > 0$;\smallskip
\item[2.] A function $\tau_s:{\cal D}_x \to [\tau_{\min}, \infty)$;\smallskip
\item[3.] A maximum allowed delay $\Delta_{\max} > 0$;\smallskip
\end{enumerate}
such that if the sequence of sampling instants ${\cal I}$ is inductively defined by
\begin{equation}
t_{k+1}=t_k+\tau_s(x(t_k))   \label{tk sequence}
\end{equation}
and if the actuation delays are such that
\begin{equation}
\Delta_k\in[0,\Delta_{\max}),\quad\forall\>k\ge0\label{constraint on Deltak}
\end{equation}
then the origin of the closed loop system~\eqref{nominal controlled nonlinear system} with control input signal $u_{{\cal I}}(t)$ is asymptotically stable.\fine
\end{problem}

\smallskip

\begin{problem} {\sl (Safety problem)}\label{Safety Problem}
Given the perturbed system~\eqref{perturbed system}, a stabilizing state feedback control law $\kappa$, and an arbitrary safe set ${\cal B}_\delta=\{x\in\R^n\mid\|x\|<\delta\}\subset{\cal D}_x$, determine $\tau_{\min}$, $\tau_s$ and $\Delta_{\max}$ as defined in Problem~\ref{Stability Problem}) such that if ${\cal I}$ is inductively defined by~\eqref{tk sequence} and if $\Delta_k$ satisfies~\eqref{constraint on Deltak}, then the closed loop system~\eqref{perturbed controlled nonlinear system} with control input signal $u_{{\cal I}}(t)$ is safe with respect to ${\cal B}_\delta$, for the time interval $[t_0,\infty)$.\fine
\end{problem}
\smallskip

In Problems \ref{Stability Problem} and \ref{Safety Problem}, the function $\tau_s$ is used to determine the next sampling instant as a function of the current measurement of the system. We require that the time interval between two sampling instants is lower bounded by a minimum sampling time $\tau_{\min} > 0$, in order to avoid undesired Zeno behaviors. We also require that the system is robust with respect to actuation delays bounded by $\Delta_{\max}$.

By choosing the next sampling instant $t_{k+1}$ as a function of the current measurement at time $t_k$, we perform sampling only when needed for guaranteeing asymptotic stability or safety. The aim is to obtain a sequence of sampling instants ${\cal I}$ with the property that the inter--sampling time $t_{k+1} - t_k$ is as large as possible, in order to reduce transmission power of the sensing and actuation data transmissions, and to reduce the CPU effort due to the computation of the control. In this paper we do not address the problem of the optimality of the solution: such a requirement can be taken into account by considering appropriate cost functions.


\section{Self Triggered Stabilizing Control}\label{sec:selfTrigStabCnt}

The results developed in this Section are based on the following assumption, analogous to the assumptions used in \cite{Tabuada 2007}.

\smallskip
\begin{assumption}\label{Assumption 1}
Assume that
\begin{enumerate}
\item[1.] $f_0 \in C^\ell({\cal D}_x \times {\cal D}_u)$, with $\ell$ a positive integer sufficiently large;
\item[2.] There exists a nonempty set ${\cal U}$ of state feedback laws $\kappa \colon{\cal D}_x\to{\cal D}_u$, such that $\kappa \in C^\ell({\cal D}_x)$ and the origin of~\eqref{nominal controlled nonlinear system} is asymptotically stable, with region of attraction a certain compact $\Omega\subset{\cal D}_x$;
\item[3.] The functions $\alpha_3,\alpha_4\in{\cal K}$ in~\eqref{Lyapunov candidate} are such that $\alpha_3^{-1},\alpha_4$ are Lipschitz.\fine
\end{enumerate}
\end{assumption}
\smallskip


The assumption of existence of a stabilizing control (i.e. non-emptiness of the set ${\cal U}$) is not restrictive, since if the nominal system cannot be stabilized using continuous time measurement and actuation, then it is clear that the nominal system cannot be stabilized using a digital control with zero--order holders. The main limitation of Assumption~\ref{Assumption 1}, and those used in \cite{Tabuada 2007}, is the Lipschitz condition on $\alpha_3^{-1}(\cdot)$ and $\alpha_4(\cdot)$. We will show how to weaken this assumption in Section~\ref{sec:selfTrigSafetyCnt}, which will be devoted to safety control. However, note that the conditions of Assumption~\ref{Assumption 1} are weaker than those used in \cite{Anta 2010} (homogeneity of the closed loop dynamics).


\medskip
\begin{theorem}\label{Th Stability for Unperturbed Systems}
Let us consider the nominal system~\eqref{nominal system}. Under Assumption~\ref{Assumption 1}, there exist a state feedback control law $\kappa$, a minimum sampling time $\tau_{\min}>0$, a function $\tau_s:{\cal D}_x \to  [\tau_{\min},\infty)$ and a maximum delay $\Delta_{\max}>0$, such that if ${\cal I}$ is inductively defined by~\eqref{tk sequence}, and $\Delta_k$ satisfies~\eqref{constraint on Deltak}, then the origin of the closed loop system~\eqref{nominal controlled nonlinear system} with control $u_{{\cal I}}(t)$ is asymptotically stable.~\fine
\end{theorem}

\proof Let us first prove the result for $\Delta_k=0$. Since ${\cal U}$ is not empty, by Assumption~\ref{Assumption 1}, we pick a state feedback control law $\kappa\in{\cal U}$. Since $f_0(x,\kappa(x))\in C^\ell({\cal D}_x)$ with $\ell>1$, there exists a Lyapunov candidate~\eqref{Lyapunov candidate}. Let us choose $r>0$ such that the ball $B_r=\{x\in\Omega\mid\norma{x}\le r\}\subset\Omega$. For $x_k\in B_r$,
\begin{equation}
\eqalign{
\dot V &= \frac{\partial V}{\partial x} f_0(x, \kappa(x_k)) = \frac{\partial V}{\partial x}f_0(x,\kappa(x)) \cr
       & \qquad +\frac{\partial V}{\partial x}\Big(f_0(x,\kappa(x_k))-f_0(x,\kappa(x))\Big)\cr
&\le -\alpha_3(\|x\|) + \alpha_4(\|x\|) \|d_h\|\cr}                                                    \label{Lyapunov inequality}
\end{equation}
where
$$
d_h = f_0(x(t), \kappa(x(t_k)))-f_0(x(t),\kappa(x(t)))
$$
is the perturbation due to the holding.

Under Assumption~\ref{Assumption 1}, there exists a $\delta_k>0$ such that $\dot x=f(x,\kappa(x_k))$ has a unique solution over $[t_k,t_k+\delta_k]$. Hence, we can expand the components $d_{h,i}$ of $d_h$ in Taylor series. Let us consider the $i^{th}$ component $d_{h,i}$, $i=1,\cdots,n$, of the $n$--dimensional vector $d_h$. One can expand each component in Taylor series with respect to $t\in[t_k,t_k+\delta_k]$, on the right of $t_k$, up to the $2^{nd}$ term, with the Lagrange remainder of the $3^{rd}$ term~\cite{Taylor 1969}
\begin{equation}
d_{h,i}=\varphi_{1,i}(x_k)(t-t_k)+\varphi_{2,i}(\bar x_i,x_k)(t-t_k)^2       \label{Taylor expansion}
\end{equation}
where
$$
\varphi_{1,i}(x_k)=\left.\d{d_{h,i}}\right|_{x(t)=x_k}\quad \varphi_{2,i}(\bar x_i, x_k)= \frac{1}{2}\left.\dd{d_{h,i}}\right|_{x(t)=\bar x_i}
$$
where $\dn{(\cdot)}$ denotes the $n$--th right derivative. According to Taylor theorem with the Lagrange remainder, there exists $\bar t_i\in[t_k,t]$, with $\bar x_i=x(\bar t_i)$, $i=1,\cdots,n$, such that the equality~\eqref{Taylor expansion} holds. Hence,
$$
\|d_h\|\le\|\varphi_{1}(x_k)\|(t-t_k)+ \|\varphi_{2}(\bar x, x_k)\|(t-t_k)^2
$$
where $\bar x\doteq(\bar x_1,\cdots,\bar x_n)$ and
$$
\eqalign{
\varphi_{1}(x_k) & \doteq \Big(\varphi_{1,1}(x_k), \cdots, \varphi_{1,n}(x_k) \Big)^T\cr
\varphi_{2}(\bar x, x_k) & \doteq \Big(\varphi_{2,1}(\bar x_1, x_k), \cdots, \varphi_{2,n}(\bar x_n, x_k) \Big)^T.\cr}
$$
Let us consider the level set $\Omega_{V(x_k)}$, and define
$$
M_1(x_k) \doteq \|\varphi_{1}(x_k)\|,\qquad  M_2(x_k)\doteq\max\limits_{\bar x\in \Omega_{V(x_k)}}\|\varphi_2(\bar x, x_k)\|.
$$
Since $f, \kappa \in C^\ell$ and $\Omega_{V(x_k)}$ is compact, then $M_1(x_k)$ is finite for any $x_k \in \Omega_{V(x_k)}$, and $M_2(x_k) \in \R^+$ exists and is finite for any $x_k \in \Omega_{V(x_k)}$.

Let us now check that there exists a time interval $[t_k,t_{k+1}<t_k+\delta_k]$ such that
\begin{equation}
\alpha_4(\|x\|)\|d_h\|\le\vartheta \alpha_3(\|x\|)           \label{Lyapunov condition}
\end{equation}
is satisfied for a fixed $\vartheta\in(0,1)$. In fact, \eqref{Lyapunov condition} is satisfied if
$$
\alpha_3^{-1}\left(\vartheta^{-1} \alpha_4(\|x\|) \left(M_1(x_k)(t-t_k)+M_2(x_k)(t-t_k)^2\right)\right) \le \|x\|.
$$
Since $\alpha_3^{-1}$ and $\alpha_4$ are Lipschitz, then equation~\eqref{Lyapunov condition} is satisfied if
$$
\vartheta^{-1} L_{\alpha_3^{-1}} L_{\alpha_4} \|x\| \left(M_1(x_k)(t-t_k)+M_2(x_k)(t-t_k)^2\right) \le  \|x\|
$$
where $ L_{\alpha_3^{-1}} > 0$ and $L_{\alpha_4} > 0$ are the Lipschitz constants respectively of $\alpha_3^{-1}$ and $\alpha_4$. The above equation directly implies that \eqref{Lyapunov condition} is satisfied if
\begin{equation}
M_1(x_k)(t-t_k)+M_2(x_k)(t-t_k)^2 \le \vartheta L^{-1}_{\alpha_3^{-1}} L^{-1}_{\alpha_4}.   \label{eqStabilityCondition}
\end{equation}
If we define
$$
\eqalign{
\tau_s(x_k)&\doteq \max\big\{t-t_k \mid \hbox{ equation \eqref{eqStabilityCondition} is satisfied}  \cr
        &\hskip3.5cm\hbox{ for each } t-t_k \in [0,\tau_s(x_k)]\big\}                                          \cr
\tau_{\min}&\doteq \min\limits_{x_k \in \Omega_{V(x_k)}}\big\{\tau_s(x_k)\big\}                       \cr}
$$
and we choose $t_{k+1}=t_k+\tau_s(x_k)$, then $\dot V(t)\le-(1-\vartheta)\alpha_3(\|x\|)$ for all $t \in [t_k, t_{k+1}]$ and for all $k\ge0$. This implies that the origin is asymptotically stable. Equation \eqref{eqStabilityCondition} is a second degree inequality in the form $a(x_k)y^2 + by \le c$, where $a(x_k),b$ are non-negative and upper bounded for each $x_k \in {{\cal D}_x}$, and $c$ is strictly positive and upper bounded. This trivially implies that $\tau_{s}(x_k)$ is strictly positive for each $x_k \in \Omega_{V(x_k)}$, and thus $\tau_{\min}$ is strictly positive as well. This completes the proof for $\Delta_k = 0$.


We now solve the problem for $\Delta_k > 0$. Following the same reasoning as above, for $t\ge t_k+\Delta_k$
$$
\eqalign{
\dot V(t) &= \frac{\partial V}{\partial x} f_0(x(t), \kappa(x_k))=\frac{\partial V}{\partial x}f_0(x,\kappa(x))\cr
       &\qquad +\frac{\partial V}{\partial x}\Big(f_0(x(t), \kappa(x(t_k+\Delta_k)))- f_0(x,\kappa(x))\Big) \cr
       &\qquad +\frac{\partial V}{\partial x}\Big(f_0(x(t), \kappa(x_k)) - f_0(x(t), \kappa(x(t_k + \Delta_k)))\Big)\cr
       &\le -\alpha_3(\|x\|) + \alpha_4(\|x\|) \|d_h\| + \alpha_4(\|x\|) \|d_{\Delta_k}\|\cr}
$$
where
$$
\eqalign{
d_h&=f_0(x(t), \kappa(x(t_k+\Delta_k)))- f_0(x,\kappa(x))\cr
d_{\Delta_k}&=f_0(x(t),\kappa(x_k))-f_0(x(t),\kappa(x(t_k + \Delta_k)))\cr}
$$
are the perturbation due to the holding and to the actuation delay. Since also the solution $x(t)$ is Lipschitz, as well as $f_0$ and $\kappa$
$$
\|d_{\Delta_k}\|\le M_3 \Delta_k,\quad M_3=L_{f_0} L_\kappa L_x
$$
where $L_{f_0}$, $L_\kappa$, $L_x$ are the Lipschitz constants of $f_0$, $\kappa$, $x$. Proceeding for $d_h$ as in the previous case, we conclude that~\eqref{Lyapunov condition} is satisfied if
\begin{equation}
M_1(x_k)(t-t_k)+M_2(x_k)(t-t_k)^2 + M_3 \Delta_k \le \vartheta L^{-1}_{\alpha_3^{-1}} L^{-1}_{\alpha_4} \label{eqStabilityConditionActuationDelay}
\end{equation}
Let $\vartheta = \vartheta_1 + \vartheta_2$, with $\vartheta_1, \vartheta_2 \in (0,1)$ and $\vartheta_1 + \vartheta_2 < 1$. We can rewrite equation \eqref{eqStabilityConditionActuationDelay} as follows
\begin{equation}
\eqalign{
M_1(x_k)&(t-t_k)+M_2(x_k)(t-t_k)^2 \cr
      & + M_3 \Delta_k \le \vartheta_1 L^{-1}_{\alpha_3^{-1}} L^{-1}_{\alpha_4} +\vartheta_2 L^{-1}_{\alpha_3^{-1}} L^{-1}_{\alpha_4}.        \cr}   \label{eqStabilityConditionActuationDelay2}
\end{equation}
Equation \eqref{eqStabilityConditionActuationDelay2} implies that the stability condition \eqref{Lyapunov condition} holds if
\begin{equation}
M_1(x_k)(t-t_k)+M_2(x_k)(t-t_k)^2 \le \vartheta_1 L^{-1}_{\alpha_3^{-1}} L^{-1}_{\alpha_4},\label{eqDeltakstabCond1}
\end{equation}
and
\begin{equation}
\Delta_k \leq (1-\vartheta_1) M_3^{-1} L^{-1}_{\alpha_3^{-1}} L^{-1}_{\alpha_4}. \label{eqDeltakstabCond2}
\end{equation}
Defining
$$
\eqalign{
\tau'_s(x_k)&\doteq \max\big\{t-t_k \mid \hbox{ equation \eqref{eqDeltakstabCond1} is satisfied}  \cr
        &\hskip3.5cm\hbox{ for each } t-t_k \in [0,\tau'_s(x_k)]\big\}                                          \cr
\tau_{\min}&\doteq \min\limits_{x_k \in \Omega_{V(x_k)}}\big\{\tau'_s(x_k)\big\}                       \cr
\Delta_{\max} &\doteq \min{\left\{(1-\vartheta_1) M_3^{-1} L^{-1}_{\alpha_3^{-1}} L^{-1}_{\alpha_4}, \tau_{\min}\right\}} \cr
\tau_s(x_k) &\doteq \tau'_s(x_k) - \Delta_{\max}}
$$
and if we choose $t_{k+1} = t_k + \tau_s(x_k)$, then $\dot V(t) \le -(1-\vartheta)\alpha_3(\|x\|)$ for all $t \in [t_k + \Delta_k, t_{k+1} + \Delta_{k+1}]$ and for all $k > 0$. This ensures the asymptotic stability of the origin. As discussed above, $\tau_{s}(x_k)$ is strictly positive for each $x_k \in \Omega_{V(x_k)}$, and thus $\tau_{\min}$ and $\Delta_{max}$ are strictly positive as well. This completes the proof.\fineTh

\begin{remark}
The choice of $\vartheta_1 \in (0,1)$ establishes an intuitive tradeoff between allowance of larger inter--sampling times $\tau_s(x_k)$ and robustness to larger actuation delays $\Delta_{\max}$. As $\vartheta_1$ decreases, $\tau_s(x_k)$ decreases and $\Delta_{\max}$ increases. This implies that we improve robustness vs actuation delays, paid by stronger sampling requirements. On the contrary, as $\vartheta_1$ increases, $\tau_s(x_k)$ increases and $\Delta_{\max}$ decreases. This implies that we achieve less restrictive sampling requirements, paid by loss of robustness to the actuation delays.
\end{remark}

\begin{remark}
When applying the self triggering rule defined in the above Theorem in a real scenario, it is necessary to solve the on-line computation of the next sampling time for each time instant $t_k$. This computation corresponds to solving a second degree equality, and is thus acceptable in an embedded system. On the contrary, the computation of $M_1(x_k)$ and $M_2(x_k)$ can be performed off-line. However, it might be difficult to compute in a closed form $M_2(x_k) = \max\limits_{\bar x \in \Omega_{V(x_k)}} \| \varphi_2(\bar x, x_k) \|$ as a function of $x_k$. In this case, we can define
$$
M_2 \doteq \max\limits_{\bar x, x_k \in \Omega_{V(x_k)}} \| \varphi_2(\bar x, x_k) \|
$$
and use it in Equation \ref{eqStabilityCondition} to compute the next sampling time.
\end{remark}

The above Remarks also apply to Theorems \ref{Th Safety for Unperturbed Systems} and \ref{Th Safety for Perturbed Systems} in following Sections.

\section{Self Triggered Safety Control}\label{sec:selfTrigSafetyCnt}

The main limitation of the results developed in Section \ref{sec:selfTrigStabCnt} is the Lipschitz continuity assumption of $\alpha_3^{-1}(\cdot)$ and $\alpha_4(\cdot)$. The following example shows that even exponentially stabilizable systems do not always satisfy this assumption.

\medskip
\begin{example}\label{exNumerical}
Consider the system $\dot x=Ax+Bu+f(x,u)=f_0(x,u)$ with
$$
f_0(x,u)=\pmatrix{-x_1+x_2+x_1^2 \cr (1+x_1)u}.
$$
Let $u=\kappa(x)=-x_2\in{\cal U}$. Consider the Lyapunov candidate $V(x)=x^T P x$, with $P$ solution of the Lyapunov equation $P A_c + A_c^T P = -Q$, with $Q=2I$ and $A_c=\pmatrix{-1 & 1 \cr 0 & -1\cr}$. Since $P=\pmatrix{2 & 1\cr 1  & 3\cr}$, then $\lambda_{\min}^P \cong 1.382$ and $\lambda_{\max}^P \cong 3.618$ denote respectively the minimum and the maximum eigenvalue of $P$. For $\|x\| \le 2/3$, the time derivative of $V$ satisfies
$$
\eqalign{
\dot V&=-\norma{x}^2_Q+2|x_1|^3+3|x_1|x_2^2 \cr
&\le -2x_1^2-2x_2^2+2 (2/3) x_1^2+3 (2/3) x_2^2 \le -\frac{1}{2}\|x\|^2\cr}
$$
thus the origin is locally exponentially stable, with $\alpha_1(\|x\|)=\lambda_{\min}^P \|x\|^2$, $\alpha_2(\|x\|)=\lambda_{\max}^P \|x\|^2$, $\alpha_3(\|x\|)=\|x\|^2/2$, and $\alpha_4(\|x\|)=\lambda_{\max}^P \|x\|$. It is clear that Assumption \ref{Assumption 1} is not satisfied, since $\alpha_3^{-1}(\cdot)$ is not Lipschitz. For this reason, we can not imply the existence of a stabilizing self triggering strategy.\fine
\end{example}
The main problem is that, if $\alpha_3^{-1}(\cdot)$ is not Lipschitz, the next sampling time $\tau_s(x_k)$ goes to zero as $x_k$ approaches the equilibrium point, and this might generate Zeno behaviors. In this Section, without the Lipschitz assumption on $\alpha_3^{-1}(\cdot)\in{\cal K}$, we will show that

\begin{enumerate}
\item[1.] For the unperturbed system \eqref{nominal system}, it is possible to keep the state arbitrarily close to the equilibrium point by applying a self triggering strategy;\smallskip

\item[2.] For the perturbed system \eqref{perturbed system}, it is possible to keep the state in a $\delta$ boundary of the equilibrium point if the disturbance norm is upper bounded by a class ${\cal K}$ function $\nu(\delta)$.
\end{enumerate}

The results developed in this section are based on the following.

\medskip
\begin{assumption}\label{assumptionSafetyUnperturbed}
Assume that $f_0 \in C^\ell({\cal D}_x \times {\cal D}_u)$, with $\ell$ a positive integer sufficiently large. Assume that there exists a nonempty set ${\cal U}$ of state feedback laws $\kappa \colon {\cal D}_x\to{\cal D}_u$, such that $\kappa \in C^\ell({\cal D}_x)$ and the origin of the system~\eqref{nominal controlled nonlinear system} is asymptotically stable.\fine
\end{assumption}


\subsection{Unperturbed Systems}

The following theorem states that, if a system is asymptotically stabilizable using a continuous time state feedback control law, then it is always possible to keep the state arbitrarily close to the equilibrium point by applying a digital self triggering strategy. Note that, in order to guarantee that the state is {\sl arbitrarily} close to the equilibrium point, we need the stabilizability assumption.


\medskip
\begin{theorem}\label{Th Safety for Unperturbed Systems}
Given the nominal system~\eqref{nominal system} and a safe set ${\cal B}_\delta$, $\delta > 0$, under Assumption~\ref{assumptionSafetyUnperturbed} there exist a state feedback control law $\kappa$, a minimum sampling time $\tau_{\min} > 0$, a function $\tau_s : {\cal D}_x \to [\tau_{\min},\infty)$ and a maximum allowed actuation delay $\Delta_{\max} > 0$, such that if ${\cal I}$ is inductively defined by $t_{k+1} = t_k + \tau_s(x(t_k))$, if $\Delta_k \in [0, \Delta_{\max})$ for each $k \ge 0$, and if the system is safe before applying the state feedback control law
$$
x(t) \in {\cal B}_\delta,\qquad \forall t \in [t_0,t_0 + \Delta_0],
$$
then the closed loop system~\eqref{nominal controlled nonlinear system} with control input signal $u_{{\cal I}}(t)$ is safe with respect to ${\cal B}_\delta$, for the time interval $[t_0, \infty)$.
\end{theorem}

\proof Using the same reasoning of Theorem \ref{Th Stability for Unperturbed Systems}, and directly considering the case $\Delta_k > 0$, we conclude that the following inequality
$$
\eqalign{
\dot V &\le -(1-\vartheta)\alpha_3(\|x\|)\!+\!\alpha_4(\|x\|)(\|d_h\|\! +\! \|d_{\Delta}\|)\! -\!\vartheta \alpha_3(\|x\|)\cr
&\le -(1-\vartheta)\alpha_3(\|x\|)\cr}
$$
holds when
$$
\alpha_4(\|x\|)\Big(M_1(x_k)(t-t_k)+M_2(x_k)(t-t_k)^2 + M_3 \Delta_k \Big) \le \vartheta \alpha_3(\|x\|)
$$
with $\vartheta \in (0,1)$, and $d_h$, $d_{\Delta}$, $M_1(x_k)$, $M_2(x_k)$, $M_3$ defined as in Theorem \ref{Th Stability for Unperturbed Systems}. The above inequality holds if
$$
\eqalign{
\|x\|\ge& \alpha_3^{-1}\bigg(\frac{\alpha_4(\delta)}{\vartheta} \Big(M_1(x_k)(t-t_k) + M_2(x_k)(t-t_k)^2\crr{-8}
&\hskip5.25cm + M_3 \Delta_k \Big) \bigg) \doteq \eta.}
$$
This implies, by~\cite {Khalil 2002}, that there exists $b := \alpha_1^{-1}(\alpha_2(\eta)) > 0$ such that $\| x(\tau)\|\le b$, $\forall \tau \in [t_k, t]$ if $x(t_k) \in {\cal B}_b$ and if the following holds
\begin{equation}
\eqalign{
&\alpha_4(\delta) \left( M_1(x_0)(t-t_k)+ M_2(x_k)(t-t_k)^2 +  M_3 \Delta_k \right)\cr
\le& \vartheta\alpha_3\Big(\alpha_2^{-1}\big(\alpha_1(\delta)\big)\Big),        \label{eqSafetyConditionUnperturbed}}
\end{equation}
where we imposed the constraint $b = \delta$. Let $\vartheta = \vartheta_1 + \vartheta_2$, with $\vartheta_1, \vartheta_2 \in (0,1)$ and $\vartheta_1 + \vartheta_2 < 1$. Equation \eqref{eqSafetyConditionUnperturbed} holds if the following inequalities hold
\begin{equation}
\eqalign{
&\alpha_4(\delta) \left( M_1(x_0)(t-t_k)+ M_2(x_k)(t-t_k)^2\right)\cr
\le &\vartheta_1\alpha_3\Big(\alpha_2^{-1}\big(\alpha_1(\delta)\big)\Big),\label{eqSafetyConditionUnperturbed1}}
\end{equation}
\begin{equation}
\alpha_4(\delta) M_3 \Delta_k \le \vartheta_2\alpha_3\Big(\alpha_2^{-1}\big(\alpha_1(\delta)\big)\Big).
\end{equation}
If we define
$$
\eqalign{
\tau'_s(x_k)&\doteq \max\big\{t-t_k \mid \hbox{ equation \eqref{eqSafetyConditionUnperturbed1} is satisfied}  \cr
        &\hbox{ for each } t-t_k \in [0,\tau'_s(x_k)]\big\}                                          \cr
\tau_{\min}&\doteq \min\limits_{x_k \in {\cal B}_{\delta}}\big\{\tau'_s(x_k)\big\}                       \cr
\Delta_{max} &\doteq \min\left\{ (1-\vartheta_1) \frac{\alpha_3\Big(\alpha_2^{-1}\big(\alpha_1(\delta)\big)\Big)}{\alpha_4(\delta) M_3}, \tau_{\min}\right\} \cr
\tau_s(x_k) &\doteq \tau'_s(x_k) - \Delta_{\max}}
$$
and we choose $t_{k+1} = t_k + \tau_s(x_k)$, then \eqref{eqSafetyConditionUnperturbed1} holds for all $t \in [t_k+ \Delta_{k}, t_{k+1}+\Delta_{k+1}]$ and for all $k \ge 0$. Since $M_1(x_k)$, $M_2(x_k)$ and $M_3$ are non-negative and upper bounded for each $x_k \in {\cal B}_{\delta}$, and since $\alpha_4(\cdot)$ and $\alpha_3\Big(\alpha_2^{-1}\big(\alpha_1(\cdot)\big)\Big)$ are ${\cal K}$ class functions, then Equation \eqref{eqSafetyConditionUnperturbed1} is a second degree inequality in the form $ay^2 + by - c \le 0$, where $a,b$ are non-negative and bounded and $c$ is strictly positive and bounded. This trivially implies that $\tau_{s}(x_k)$ is strictly positive for each $x_k \in {\cal B}_\delta$, and thus $\tau_{\min}$ is strictly positive as well. It is also clear that $\Delta_{\max}$ is strictly positive. This completes the proof.~\quad~\fineTh
\medskip

%


\subsection{Perturbed Systems}

A generic system~\eqref{perturbed system}, subject to disturbances and parameter variations, can be seen as the nominal system~\eqref{nominal system}, perturbed by the term
\begin{equation}
g(x,u,\mu,d) = f(x,u,\mu,d)-f_0(x,u) \doteq d_g.         \label{eqPerturbation}
\end{equation}
Hence, \eqref{perturbed system} can be rewritten as follows
\begin{equation}
\dot x=f_0(x,u)+g(x,u,\mu,d).
\end{equation}
\medskip
\begin{definition}\label{admissible perturbation}
Under Assumption~\ref{assumptionSafetyUnperturbed}, and given the perturbed system~\eqref{perturbed system} and a safe set ${\cal B}_\delta$, $\delta > 0$, we define the perturbation~\eqref{eqPerturbation} $\delta$--admissible if there exists a state feedback control law $\kappa\in{\cal U}$ and a constant $\vartheta_g\in(0,1)$ such that the function $g(x,\kappa(x_0),\mu,d)$ satisfies
\begin{equation}
\max_{\phantom{.}\atop{x,x_k\in{\cal B}_\delta\atop{d\in{\cal D}_d\atop \mu\in{\cal D}_\mu}}}\!\!\!\!\|g(x,\kappa(x_k),\mu,d)\|\!\le\! \vartheta_g\frac{\alpha_3\Big(\alpha_2^{-1}\!\big(\alpha_1(\delta)\big)\Big)}{\alpha_4(\delta)}\! \doteq\! \nu(\delta) \quad \label{admissibility}
\end{equation}
with $\alpha_1,\alpha_2,\alpha_3,\alpha_4$ as in~\eqref{Lyapunov candidate}.\fine
\end{definition}
The $\delta$--admissible perturbations are those for which the safety problem with respect to a ball ${\cal B}_\delta$ can be solved using continuous time measurement and actuation. If a perturbation is not $\delta$--admissible, safety with respect to ${\cal B}_\delta$ is clearly not achievable using sampled measurements and actuations. Note that in condition~\eqref{admissible perturbation} the expression of $\nu(\delta)$ can be explicitly computed.

The following theorem states that, if a system is asymptotically stabilizable using a continuous time state feedback control law and the perturbation is $\delta$--admissible, then it is possible to keep the state in a boundary ${\cal B}_\delta$ of the equilibrium point by applying a digital self triggering strategy.

\medskip
\begin{theorem}\label{Th Safety for Perturbed Systems}
Given the perturbed system~\eqref{perturbed system} and a safe set ${\cal B}_\delta$, $\delta > 0$, under Assumption~\ref{assumptionSafetyUnperturbed} and for any $\delta$--admissible perturbation~\eqref{eqPerturbation}, there exist a state feedback control law $\kappa$, a minimum sampling time $\tau_{\min} > 0$, a function $\tau_s:{\cal D}_x \to [\tau_{\min},\infty)$ and a maximum allowed actuation delay $\Delta_{\max} > 0$, such that if ${\cal I}$ is inductively defined by $t_{k+1} = t_k + \tau_s(x(t_k))$, if $\Delta_k \in [0, \Delta_{\max})$ for each $k \ge 0$, and if the system is safe before applying the state feedback control law
$$
x(t) \in {\cal B}_\delta,\qquad \forall t \in [t_0,t_0 + \Delta_0]
$$
then the closed loop system~\eqref{perturbed controlled nonlinear system} with control input signal $u_{{\cal I}}(t)$ is safe with respect to ${\cal B}_\delta$, for the time interval $[t_0 + \Delta_0,\infty)$.~\quad~\fine
\end{theorem}

\proof Using the same reasoning of Theorem \ref{Th Safety for Unperturbed Systems} , we conclude that the following inequality
$$
\eqalign{
\dot V &\le -(1-\vartheta)\alpha_3(\|x\|) -\vartheta \alpha_3(\|x\|) \cr
&\qquad\qquad +\alpha_4(\|x\|)(\|d_h\| + \|d_{\Delta}\| + \|d_{g}\|) \cr
&\le -(1-\vartheta)\alpha_3(\|x\|)\cr}
$$
holds when
$$
\eqalign{
\alpha_4(\|x\|)&\Big(M_1(x_k)(t-t_k)+M_2(x_k)(t-t_k)^2 + M_3 \Delta_k + \|d_{g}\|\Big)\cr
&\le \vartheta \alpha_3(\|x\|)}
$$
with $\vartheta \in (0,1)$, and $d_h$, $d_{\Delta}$, $M_1(x_k)$, $M_2(x_k)$, $M_3$ defined as in Theorem \ref{Th Stability for Unperturbed Systems}. The above inequality holds is
$$
\eqalign{
\|x\| \ge& \alpha_3^{-1}\bigg(\frac{\alpha_4(\delta)}{\vartheta} \Big(M_1(x_k)(t-t_k) + M_2(x_k)(t-t_k)^2 \crr{-4}
&+ M_3 \Delta_k + \|d_{g}\| \Big) \bigg) \doteq \eta.}
$$
This implies, by~\cite {Khalil 2002}, that there exists $b := \alpha_1^{-1}(\alpha_2(\eta)) > 0$ such that $\| x(\tau)\|\le b$, $\forall \tau \in [t_k, t]$ if $x(t_k) \in {\cal B}_b$ and if the following holds
\begin{equation}
\eqalign{
\alpha_4(\delta)& \Big( M_1(x_0)(t-t_k) + M_2(x_k)(t-t_k)^2 + M_3 \Delta_k \Big) + \|d_{g}\|\cr
&\le\vartheta\alpha_3\Big(\alpha_2^{-1}\big(\alpha_1(\delta)\big)\Big).\label{eqSafetyConditionPerturbed}}
\end{equation}
where we imposed the constraint $b = \delta$. Let $\vartheta = \vartheta_1 + \vartheta_2 + \vartheta_g$, with $\vartheta_1, \vartheta_2, \vartheta_g \in (0,1)$ and $\vartheta_1 + \vartheta_2 < 1 - \vartheta_g$. Since the perturbation is assumed $\delta$--admissible, then
$$
\|d_g\| \le \vartheta_g\frac{\alpha_3\Big(\alpha_2^{-1}\big(\alpha_1(\delta)\big)\Big)}{\alpha_4(\delta)}.
$$
Thus, Equation \eqref{eqSafetyConditionPerturbed} holds if the following inequalities hold
\begin{equation}
\eqalign{
\alpha_4(\delta)& \Big( M_1(x_0)(t-t_k) + M_2(x_k)(t-t_k)^2 \Big)\cr
&\le \vartheta_1\alpha_3\Big(\alpha_2^{-1}\big(\alpha_1(\delta)\big)\Big)\cr
\alpha_4(\delta)&M_3 \Delta_k \le \vartheta_2\alpha_3\Big(\alpha_2^{-1}\big(\alpha_1(\delta)\big)\Big).
}\label{eqSafetyConditionUnperturbed3}
\end{equation}
If we define
$$
\eqalign{
\tau'_s(x_k)&\doteq \max\big\{t-t_k \mid \hbox{ equation \eqref{eqSafetyConditionUnperturbed3} is satisfied}  \cr
        &\hbox{ for each } t-t_k \in [0,\tau'_s(x_k)]\big\}                                          \cr
\tau_{\min}&\doteq \min\limits_{x_k \in {\cal B}_{\delta}}\big\{\tau'_s(x_k)\big\}                       \cr
\Delta_{max} &\doteq \min\left\{ (1-\vartheta_1-\vartheta_g) \frac{\alpha_3\Big(\alpha_2^{-1}\big(\alpha_1(\delta)\big)\Big)}{\alpha_4(\delta) M_3}, \tau_{\min}\right\} \cr
\tau_s(x_k) &\doteq \tau'_s(x_k) - \Delta_{\max}}
$$
and we choose $t_{k+1} = t_k + \tau_s(x_k)$, then \eqref{eqSafetyConditionPerturbed} holds for all $t \in [t_k+ \Delta_{k}, t_{k+1}+\Delta_{k+1}]$ and for all $k \ge 0$. As discussed in Theorem \ref{Th Stability for Unperturbed Systems}, $\tau_{s}(x_k)$ is strictly positive for each $x_k \in {\cal B}_\delta$, and $\tau_{\min}$, $\Delta_{\max}$ are strictly positive as well. This completes the proof.\fineTh

As discussed in Section \ref{sec:selfTrigStabCnt}, the choice of $\vartheta_1$, $\vartheta_2$ and $\vartheta_g$ establishes in an intuitive way the tradeoff between allowance of larger inter--sampling times ($\vartheta_1$), and robustness to larger actuation delays ($\vartheta_2$) and perturbations ($\vartheta_g$).

Theorems~\ref{Th Safety for Unperturbed Systems} and~\ref{Th Safety for Perturbed Systems} prove the existence of a self triggering strategy characterized by the time sequence ${\cal I}=\{t_k\}_{k\ge 0}$, with  $t_k \ge \tau_{\min} > 0$ for each $k\ge 0$, such that the closed loop system satisfies a given safety specification. Moreover, they provide a formula to explicitly compute the next sampling time $t_{k+1}$ as a function of the state $x(t_k)$ at time $t_k$.

Although the simulation results illustrated in Section \ref{sec:simulations} show strong benefits of the proposed self triggering strategy with respect to controllers based on constant sampling, the sequence ${\cal I}$ might be conservative, in the sense that longer sampling times might be determined, because of the approximations used in the proof. A trivial way to obtain a less conservative sequence ${\cal I}$ without introducing more restricting assumptions is the use of Taylor expansions of order higher than 2.

\Omesso

\section{Conditions for the Existence of Less Conservative Sampling Time Sequences}

Theorems~\ref{Th Safety for Unperturbed Systems} and~\ref{Th Safety for Perturbed Systems} prove the existence of a self triggering strategy characterized by the time sequence ${\cal I}=\{t_k\}_{k\ge 0}$, with  $t_k \ge \tau_{\min} > 0$ for each $k\ge 0$, such that the closed loop system satisfies a given safety specification. Moreover, it provides a formula to explicitly compute the next sampling time $t_{k+1}$ as a function of the state $x(t_k)$ at time $t_k$.

Although the simulation results illustrated in Section \ref{sec:simulations} show strong benefits of the proposed self triggering strategy with respect to controllers based on constant sampling, the sequence ${\cal I}$ might be conservative, in the sense that longer sampling times might be determined, because of the approximations used in the proof. A trivial way to obtain a less conservative sequence ${\cal I}$ without introducing more restricting assumptions is the use of Taylor expansions of order higher than 2 to obtain. The remainder of this section is devoted to the analysis of a less conservative sequence ${\cal I}$, on the basis of the following assumption.

\begin{assumption}\label{A. V is class ell}
The functions $V(x)$ and $\alpha_3(\alpha_2^{-1}(V(x)))$ are of class $C^\ell({\cal D}_x)$, where $V(x)$ is the Lyapunov function satisfying~\eqref{Lyapunov candidate}.\fine
\end{assumption}


For the sake of conciseness, we study only the case without actuation delay and perturbations: actuation delays and disturbances can be then considered by proceeding as in Theorems \ref{Th Safety for Unperturbed Systems} and \ref{Th Safety for Perturbed Systems}.

Let us consider inequality~\eqref{Lyapunov inequality}. Since $-\alpha_3(\|x\|) \le\allowbreak - \alpha_3(\alpha_2^{-1}(V))$, one gets
$$
\dot V \le - \alpha_3(\alpha_2^{-1}(V)) + \frac{\partial V}{\partial x}\Big(f_0(x, \kappa(x_k)) - f_0(x,\kappa(x))\Big).
$$
Since the product of $C^\ell({\cal D}_x)$ functions is a $C^\ell({\cal D}_x)$ function, then Assumption \ref{A. V is class ell} implies that the right side of the above inequality is a scalar function of class $C^\ell({\cal D}_x)$, and it is thus possible to expand in Taylor series with respect to the time $t$, around $t_k$, up to the $(\ell-1)^{th}$ term and with Lagrange remainder of the $\ell^{th}$ term:
\begin{equation}
\dot V \le - \alpha_3(\alpha_2^{-1}(V(x(t_k)))) + \phi(t-t_k,x(t_k))
\end{equation}
with
$$
\phi(t-t_k,x(t_k)) \doteq \sum\limits_{i=1}^{\ell} \frac{\varphi_i(x(t_k))}{i!}(t-t_k)^i.
$$
By solving the differential equation above, we obtain
\begin{equation}
\eqalign{
V(x(t)) &\le V(x(t_k))-\alpha_3(\alpha_2^{-1}(V(x(t_k))))(t-t_k) \crr{4}
                    & \qquad \qquad \qquad \qquad + \varphi(t-t_k,x(t_k))   \cr}
\end{equation}
where
$$
\phi(t-t_k,x(t_k)) \doteq \sum\limits_{i=1}^{\ell} \frac{\varphi_i(x(t_k))}{(i+1)!}(t-t_k)^{i+1}.
$$
Since $\|x(t)\|\le\alpha_1^{-1}(V(x(t)))$, then
$$
\eqalign{
\|x(t)\| \le\alpha_1^{-1}\Big[ V(x(t_k))&-\alpha_3(\alpha_2^{-1}(V(x(t_k))))(t-t_k)\cr
&\qquad\qquad\qquad + \varphi(t-t_k,x(t_k)) \Big].          }
$$
Imposing the safety condition
\begin{equation}
\eqalign{
V(x(t_k)) &- \alpha_3(\alpha_2^{-1}(V(x(t_k))))(t-t_k)\crr{4}
&\qquad\qquad \qquad + \varphi(t-t_k,x(t_k))\le\alpha_1(\delta)          \cr}             \label{equationTrigTime}
\end{equation}
for all $t\in[t_k,t_{k+1}]$, with $t_{k+1}$ the next sampling time computed such that~\eqref{equationTrigTime} is satisfied
$$
t_{k+1}=\max\{t\mid \eqref{equationTrigTime} \hbox{is safisfied}\}
$$
\begin{equation}
\tau_s(x_k) = \min \{t > 0 :{\rm Equation \eqref{equationTrigTime} is NOT safisfied}\}.
\end{equation}
Notice that $\tau_s(x_k)$ is strictly greater than zero if and only if the following holds:
\begin{equation}
\forall x(t_k) \in \cal B_{\delta},V(x(t_k)) - \alpha_3(\alpha_2^{-1}(V(x(t_k)))) < \alpha_1(\delta).
\end{equation}

%

\EndOmesso


\section{An Example of Application of the Digital Self Triggered Robust Control}\label{sec:simulations}

Consider the system defined in Example~\ref{exNumerical}. As already shown, we can not imply the existence of a stabilizing self triggering strategy. However, since Assumption \ref{assumptionSafetyUnperturbed} holds, Theorem \ref{Th Safety for Unperturbed Systems} implies the existence of a self triggering strategy that guarantees safety for an arbitrary small neighborhood of the equilibrium point. In particular, since the origin of the system is locally exponentially stabilizable for $\|x\| \le 2/3$, we define the safe set ${\cal B}_\delta$ with $\delta=10^{-4} < 2/3$. We performed simulations using Matlab, with initial condition $x_0 = (10^{-5},10^{-5})^T \in {\cal B}_\delta$.

\begin{figure}[ht]
\centerline{\includegraphics[scale=0.45]{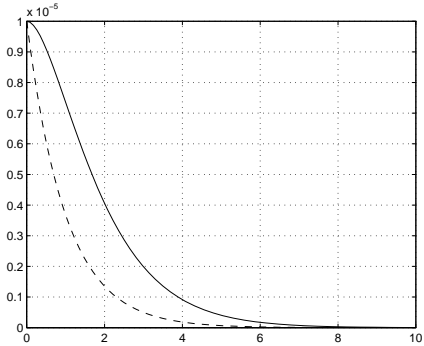}}
\caption{Continuous control: $x_1$ (solid), $x_2$ (dashed) vs time}\label{figSimResContCtrl}
\end{figure}

In Figure \ref{figSimResContCtrl}, the closed loop behavior is illustrated when a continuous time control law is used. The closed loop system is asymptotically stable.

\begin{figure}[ht]
\centerline{\includegraphics[scale=0.6]{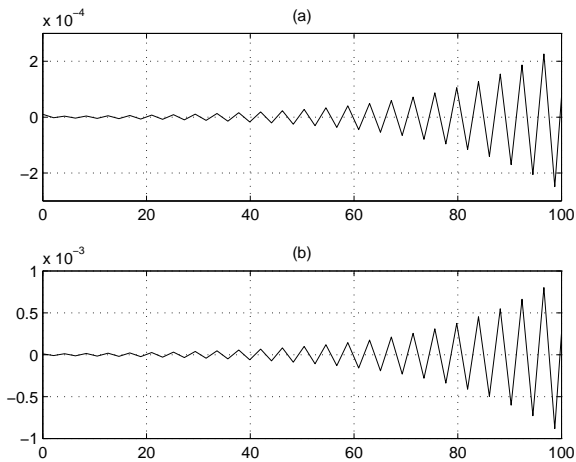}}
\caption{Digital control with constant sampling time of~2.1~s: (a)~$x_1$; (b)~$x_2$ vs time}\label{figSimResConstCtrl}
\end{figure}

In Figure \ref{figSimResConstCtrl}, the closed loop behavior is illustrated when a discrete time control law with constant sampling time of~2.1~s is used. The closed loop system is unstable.

\begin{figure}[ht]
\centerline{\includegraphics[scale=0.6]{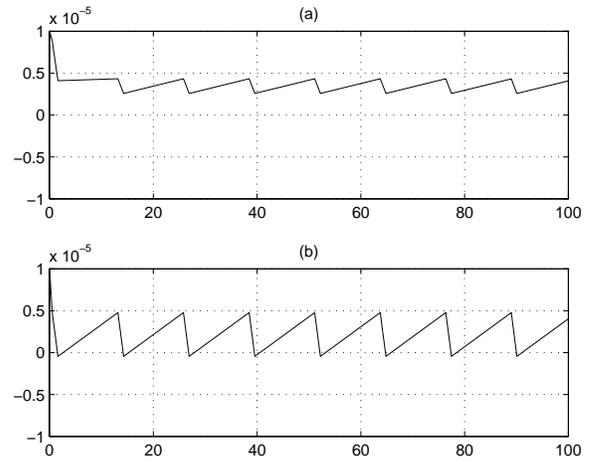}}
\caption{Self triggering control with $\vartheta=0.99$ and $\Delta_k=0$~ms: (a)~$x_1$; (b)~$x_2$ vs time}\label{figSimulationResultsNoDelay}
\end{figure}

\begin{figure}[ht]
\centerline{\includegraphics[scale=0.6]{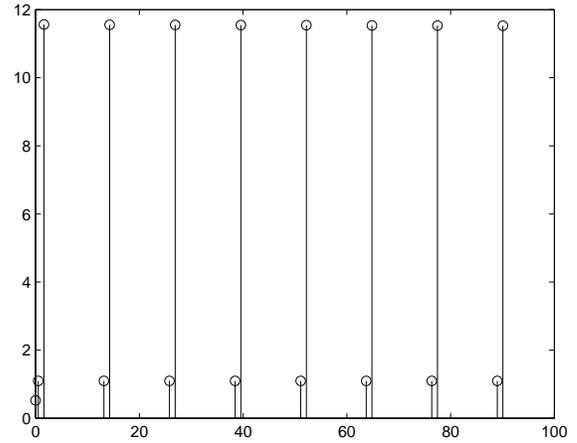}}
\caption{Sequence of sampling instants ${\cal I} = \{t_k\}_{k \ge 0}$ [s] with $\vartheta=0.99$ and $\Delta_k=0$~ms.}\label{figSamplingTimesNoDelay}
\end{figure}

In Figure \ref{figSimulationResultsNoDelay}, the closed loop behavior is illustrated when the proposed self triggering control algorithm is used, with $\vartheta = 0.99$ and with no actuation delay, namely $\Delta_k = 0$ for each $k \ge 0$ . The closed loop system is not asymptotically stable, but is safe with respect to ${\cal B}_{\delta}$ for the time interval $[t_0, \infty)$. It is interesting to remark that the average sampling time is ~6.2~s, i.e. more than $295 \%$ longer than the constant sampling time of 2.1~s that yields an unstable control loop. Thus, using the proposed self triggering control algorithm, we achieve safety reducing of more than $295 \%$ the battery energy consumption, with respect to an unstable control strategy with constant sampling. However, since we have chosen $\vartheta = 0.99$, we can only guarantee robustness with respect to actuation delays bounded by $\Delta_{\max} = 0.17$~ms.

\begin{figure}[ht]
\centerline{\includegraphics[scale=0.6]{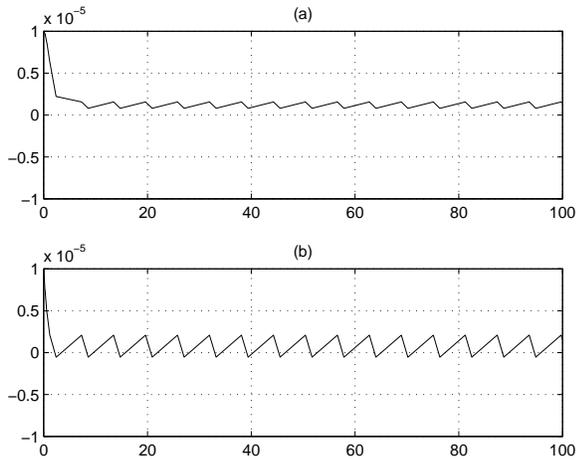}}
\caption{Self triggering control with $\vartheta = 0.5$ and $\Delta_k = 9$~ms: (a)~$x_1$; (b)~$x_2$ vs time}\label{figSimulationResultsDelay}
\end{figure}

\begin{figure}[ht]
\centerline{\includegraphics[scale=0.6]{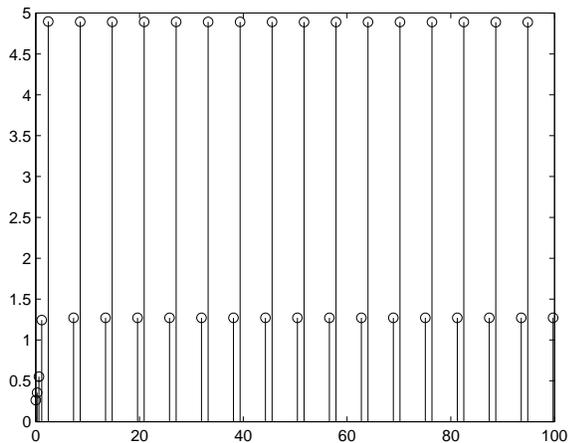}}
\caption{Sequence of sampling instants ${\cal I} = \{t_k\}_{k \ge 0}$ [s] with $\vartheta = 0.5$ and $\Delta_k = 9$~ms} \label{figSamplingTimesDelay}
\end{figure}

In Figure \ref{figSimulationResultsDelay}, the closed loop behavior is illustrated when the proposed self triggering control algorithm law is used, with $\vartheta = 0.5$ and with actuation delay $\Delta_k = \Delta_{\max} = 9$~ms for each $k \ge 0$. The closed loop system is not asymptotically stable, but is still safe with respect to ${\cal B}_{\delta}$ for the time interval $[t_0, \infty)$. However, since we have chosen $\vartheta = 0.5$ in order to be robust with respect to actuation delays, the average sampling time 3~s is more conservative with respect to the case $\vartheta = 0.99$. However, the average sampling time is almost $50 \%$ longer than the constant sampling time of 2.1~s, that yields an unstable control loop. Thus, using the proposed self triggering control algorithm, we achieve safety reducing of almost $50 \%$ the battery energy consumption, with respect to an unstable control strategy with constant sampling, while guaranteeing robustness with respect to actuation delays bounded by $\Delta_{\max} = 9$~ms.

                                                         \Omesso

\section{The Digital Self Triggered Robust Control Applied to the Vehicle Control}\label{sec:automotiveCaseStudy}

The complex dynamics of a vehicle can be approximated, under suitable hypotheses, with the so--called bicycle model. Despite its simplicity, this model captures the vehicle major characteristics of interest and is very used in the literature and in the industry to approach the design of a handling and stability controller~\cite{Bianchi 2010}
\begin{equation}
\eqalign{
\dot v_y&=-v_x\omega_z+\frac{\mu}{m}\Big(F_{y,f}(\delta_d,v_y,\omega_z,\delta_c)+F_{y,r}(v_y,\omega_z)\Big)      \crr{4}
\dot\omega_z&=\frac{\mu}{J}\Big(F_{y,f}(\delta_d,v_y,\omega_z,\delta_c)l_f-F_{y,r}(v_y,\omega_z)l_r\Big)+\frac{1}{J}M_z      \cr}\label{modello 1}
\end{equation}
where $v_y$ is the lateral velocity and $\omega_z$ is the the yaw velocity, while $m$ is the vehicle mass, $J$ is its inertia momentum, $v_x$ is the (constant) longitudinal velocity, $l_f,l_r$ are the distances from the center of gravity to the front and rear axles, $\mu$ is the tire--road friction coefficient, and $\delta_d$ is the driver angle. In~\eqref{modello 1}, some actuators are considered, in order to design an active control of the vehicle: the Active Front Steer $\delta_c$ (AFS), imposing an incremental steer angle on in the front steering angle, and the Rear Torque Vectoring $M_z$ (RTV), which distributes the torque in the rear axle. This active control is exerted by the front, rear tire lateral forces
$$
\eqalign{
F_{y,f}&=C_{y,f}\sin(A_{y,f}\arctan(B_{y,f}\alpha_f))      \crr{8}
F_{y,r}&=C_{y,r}\sin(A_{y,r}\arctan(B_{y,r}\alpha_r))      \crr{0}
\alpha_f&=\delta_d+\delta_c-\frac{v_y+l_f\omega_z}{v_x}    \crr{-4}
\alpha_r&=-\frac{v_y-l_r\omega_z}{v_x}                     \cr}
$$
where an approximated expression of the so--called Pacejka formula has been considered~\cite{Pacejka 2005}, making use of positive experimental parameters $A_{y,f}$, $B_{y,f}$, $C_{y,f}$, $A_{y,r}$, $B_{y,r}$, $C_{y,r}$. Note that $F_{y,f}$ has a saturation value for a certain value $\alpha_{f,{\rm sat}}$ of the front slip angle $\alpha_f$.

Under the assumption that $|\alpha_f|\le\alpha_{f,{\rm sat}}$, i.e. when the front tire force does not saturate, and that $\mu>0$, i.e. when the friction coefficient is non--zero, we want here to design a controller such that
$$
\lim\limits_{t\to\infty} v_y=v_{y,\reference},\quad \lim\limits_{t\to\infty} \omega_z=\omega_{z,\reference}
$$
for all initial conditions $v_y(0)$, $\omega_z(0)$, with $v_{y,\reference}(t)$, $\omega_{z,\reference}(t)$ some bounded reference functions, with bounded derivatives. These functions describe a desired behavior of the vehicle, and can be generated for instance by the reference system
$$
\eqalign{
\dot v_{y,\reference}     &=-v_x\omega_{z,\reference}+\frac{\mu}{m}(F_{y,f,\reference}+F_{y,r,\reference})\crr{4}
\dot\omega_{z,\reference} &= \frac{\mu}{J}(F_{y,f,\reference}l_f-F_{y,r,\reference}l_r)\cr}
$$
where
$$
\eqalign{
F_{y,f,\reference}&= C_{y,f,\reference} \sin(A_{y,f,\reference}\arctan(B_{y,f,\reference}\alpha_{f,\reference}))        \crr{8}
F_{y,r,\reference}&= C_{y,r,\reference} \sin(A_{y,r,\reference}\arctan(B_{y,r,\reference}\alpha_{r,\reference}))        \crr{0}
\alpha_{f,\reference}&= \delta_d-\frac{v_{y,\reference}+l_f\omega_{z,\reference}}{v_x} \crr{-4}
\alpha_{r,\reference}&=-\frac{v_{y,\reference} - l_r \omega_{z,\reference}}{v_x}       \cr}
$$
with $A_{y,f,\reference}$, $B_{y,f,\reference}$, $C_{y,f,\reference}$, $A_{y,r,\reference}$, $B_{y,r,\reference}$, $C_{y,r,\reference}$ parameters ensuring an appropriate behavior to the reference dynamics.

Setting $e_{v_y}=v_y-v_{y,\reference}$, $e_{\omega_z}=\omega_z-\omega_{z,\reference}$, the error dynamics are
$$
\eqalign{
\dot e_{v_y}&=-v_x e_{\omega_z}+\frac{\mu}{m}(F_{y,f}+F_{y,r}-F_{y,f,\reference}-F_{y,r,\reference})\crr{0}
\dot e_{\omega_z}&=\frac{\mu}{J}\Big((F_{y,f}-F_{y,f,\reference})l_f-(F_{y,r}-F_{y,r,\reference})l_r\Big)+\frac{1}{J}M_z. \cr}
$$
Hence, considering the Lyapunov candidate
$$
\alpha_1(\|e\|) \le V=\frac{1}{2\mu}(m e_{v_y}^2+J e_{\omega_z}^2)\le \alpha_2(\|e\|)
$$
with
$$
\eqalign{
\alpha_1(\|e\|)&=\min\bigg\{\frac{m}{2\mu},\frac{J}{2\mu}\bigg\}\|e\|^2\crr{0}
\alpha_2(\|e\|)&=\max\bigg\{\frac{m}{2\mu},\frac{J}{2\mu}\bigg\}\|e\|^2\crr{0}
\left\|\frac{\partial V(e)}{\partial e}\right\|\le\alpha_4(\|e\|)&=\max\bigg\{\frac{m}{\mu},\frac{J}{\mu}\bigg\}\|e\|\cr}
$$
one gets
$$
\eqalign{
\dot V&=-\frac{m}{\mu} v_x e_{v_y} e_{\omega_z}+e_{v_y}(F_{y,f}+F_{y,r}-F_{y,f,\reference}-F_{y,r,\reference})\crr{4}
      &\qquad+ e_{\omega_z}\Big((F_{y,f}-F_{y,f,\reference})l_f-(F_{y,r}-F_{y,r,\reference})l_r\Big)\crr{4}
      &\qquad +\frac{1}{\mu} e_{\omega_z}M_z.\cr}
$$
The following state feedback
$$
\eqalign{
\delta_c &= -\delta_d + \frac{v_y + l_f \omega_z}{v_x} + B_{y,f}^{-1}\tan(A_{y,f}^{-1}\arcsin(C_{y,f}^{-1}\gamma))\crr{0}
M_z &= -k_2 e_{\omega_z}-\mu\Big((F_{y,f}-F_{y,f,\reference})l_f-(F_{y,r}-F_{y,r,\reference})l_r\Big)             \cr}
$$
with $k_1,k_2>0$ and
$$
\gamma=-k_1 e_{v_y}-F_{y,r}+F_{y,f,\reference}+F_{y,r,\reference}+\frac{m}{\mu} v_x e_{\omega_z}
$$
ensures that
$$
\dot V=-k_1 e_{v_y}^2 -k_2\frac{1}{\mu} e_{\omega_z}^2\le -\alpha_3(\|e\|)=-\min\bigg\{k_1,\frac{k_2}{\mu}\bigg\} \|e\|^2
$$
i.e. the exponential tracking of the references.

We have implemented this state feedback using the proposed digital self triggered algorithm. The vehicle parameters are given in Table~1.
$$
\eqaligntwo{
v_x &= 30\ {\rm m/s}           &   A_{y,f} &= 7.2\cr
m &= 1880\ {\rm Kg}            &   B_{y,f} &= 1.81\cr
J &= 2830\ {\rm Kg\,m}^2       &   C_{y,f} &= 8854\ {\rm N}\cr
l_f &= 1.38\ {\rm m}           &   A_{y,r} &= 11\cr
l_r &= 1.53\ {\rm m}           &   B_{y,r} &= 1.68\cr
\mu &= 0.9                     &   C_{y,r} &= 8394\ {\rm N}\crr{4}}
$$
\centerline{Table 1 -- Vehicle parameters}

As already shown, we can not imply the existence of a stabilizing self triggering strategy. However, since Assumption~\ref{assumptionSafetyUnperturbed} holds, Theorem \ref{Th Safety for Unperturbed Systems} implies the existence of a self triggering strategy that guarantees safety for an arbitrary small neighborhood of the equilibrium point. We consider {\bf zero initial conditions} for $v_y$, $\omega_z$, $v_{y,\reference}$ and $\omega_{z,\reference}$.

\begin{figure}[ht]
\centerline{\includegraphics[scale=0.6]{.figures/automotivefig3.eps}}
\caption{Digital control with constant sampling time of~0.18~ms: (a)~$e_{v_y}$ [m/s vs s]; (b)~$e_{\omega_z}$ [ras/s vs s]}\label{Figure1}
\end{figure}

\begin{figure}[ht]
\noindent{\bf Togliere la scritta dalla figura}\vskip-2pt
$\longrightarrow$\vskip-2pt

\centerline{\includegraphics[scale=0.6]{.figures/automotivefig4.eps}}
\caption{Digital control with constant sampling time of~0.18~ms: Car trajectory in the inertial frame (solid) and reference (dashed)}\label{Figure2}
\end{figure}

In Fig.~\ref{Figure1}, the closed loop behavior of the error dynamics is shown when a discrete time control law with constant sampling time of~0.18~ms is used and no actuation delays are considered. The closed loop system is clearly unstable. In Fig.~\ref{Figure2} the car trajectory in the inertial frame the is shown.

\begin{figure}[ht]
\centerline{\includegraphics[scale=0.6]{.figures/automotivefig5.eps}}
\caption{Self triggering control with $\vartheta=0.99$ and $\Delta_k=0$: (a)~$e_{v_y}$ [m/s vs s]; (b)~$e_{\omega_z}$ [ras/s vs s]}\label{Figure3}

\noindent{\bf (a) Ricontrollare}

\noindent{\bf Usare {\tt axis}}
\end{figure}

\begin{figure}[ht]
\noindent{\bf Togliere la scritta dalla figura}\vskip-2pt
$\longrightarrow$\vskip-2pt

\centerline{\includegraphics[scale=0.6]{.figures/automotivefig7.eps}}
\caption{Self triggering control: Car trajectory in the inertial frame (solid) and reference (dashed)}\label{Figure4}

\noindent{\bf Andamento impossibile. Rifare}
\noindent{\bf Usare {\tt axis}}
\end{figure}

When the proposed self triggering control algorithm is used, with $\vartheta = 0.99$, and no actuation delay is present, i.e. $\Delta_k = 0$, the closed loop system is not asymptotically stable, but is safe with respect to ${\cal B}_{0.4}$ for the time interval $[t_0, \infty)$, see Figs.~\ref{Figure3}, \ref{figure4}. It is interesting to remark that the average sampling time is~0.32~ms, i.e. more than $175\%$ longer than the constant sampling  time of 0.18~ms that yields unstable error dynamics. Thus, using the proposed self triggering control algorithm, we achieve safety reducing of more than $175\%$ the battery energy consumption, with respect to an unstable control strategy with constant sampling.

Similar results can be obtained in the presence of $\delta$--admissible perturbations or delays. Alternatively, one can that care of more severe perturbations or delays at the expense of bigger $\delta$, or of higher values of the sampling time.

                                                         \EndOmesso

\section{Conclusions}\label{sec:conclusions}

We have developed novel results on self triggering control for the asymptotic stability of unperturbed nonlinear systems, affected by bounded actuation delays, and for the safety problem, for nonlinear systems perturbed by norm--bounded parameter uncertainties and disturbances, and affected by bounded actuation delays. We have provided a methodology for the computation of the next execution time in both cases. We have showed on a simple case study that the proposed results provide strong benefits in terms of energy consumption, with respect to digital controls based on constant samplings, by reducing the average sampling times. As a next step of this research line, we aim to tackle more complex case studies, and obtain results for less conservative sampling time sequences.

\end{document}